\documentclass[9pt]{article}
\usepackage{amssymb, amsmath, amsthm, amstext, amscd, amsgen, amsfonts}
\usepackage{latexsym, wasysym, graphics, mathrsfs}
\usepackage{bm, bbm}
\usepackage{euscript, enumerate, enumitem, calc, verbatim}
\usepackage{microtype, colonequals} % I have added this package from the suggestion by Bjorn Poonen homepage
\usepackage{tikz}
\usetikzlibrary{cd}
\usetikzlibrary{decorations.pathmorphing}
\usepackage[all]{xy}
\usepackage{color}
\usepackage[colorlinks, pagebackref]{hyperref}
\hypersetup{
colorlinks=true,
citecolor=blue,
linkcolor=blue,
urlcolor=blue}
\usepackage{url} %the modfications that follows is from mathstack exchnage. Forgot to add the link here.
\makeatletter
\def\Url@twoslashes{\mathchar`\/\@ifnextchar/{\kern-.2em}{}}
\g@addto@macro\UrlSpecials{\do\/{\Url@twoslashes}}
\makeatother

\newtheorem{theorem}{Theorem}[section]
\newtheorem{lemma}[theorem]{Lemma}

\newtheorem{sub-lemma}[theorem]{Sub-lemma}

\theoremstyle{definition}

\newcommand{\suchthat}{\;\ifnum\currentgrouptype=16 \middle\fi|\;}
\newcommand{\mycomment}[1]{}

\AtEndDocument{

\bigskip{\footnotesize%
Ankit Rai \par
\textsc{Department of Mathematics, University at Buffalo, Buffalo, NY, 14260.}\par
\textit{E-mail address} : \texttt{arai22@buffalo.edu}.}}

\DeclareMathOperator*{\Spec}{\mathrm{Spec}}
\DeclareMathOperator*{\im}{\mathrm{im}}

\title{A note on varieties of non-negative Kodaira dimension with polarized self maps}
\author{Ankit Rai}

\begin{document}
\maketitle

\begin{abstract}
In this note we prove that a smooth projective variety (defined over a field $k$) of non-negative Kodaira dimension that has a $k$-rational point and a polarized self map must be a finite free quotient of an abelian variety.
\end{abstract}

\section{Introduction}\label{sec1}

	Let $X$ be a smooth projective variety over an algebraically closed field $k$. A tuple $(X, \phi, L)$ is said to be a polarized variety or a variety $X$ equipped with a polarized endomorphism $(\phi, L)$ if there exists an ample line bundle $L$ on $X$ and a self map $\phi$ such that $\phi^*(L) \cong L^{\otimes d}$ for some $d>1$.
	The study of periodic (and more generally preperiodic) points of such a dynamical system has gathered much interest in recent years resulting in progress towards Bogomolov conjecture, dynamical Mordell--Lang conjectures, etc.
	On the other hand, study of geometric properties of such varieties are also of interest and has been pursued in several article including those by Nakayama, Nakayama--Zhang, Meng--Zhang, Yoshikawa and others. The study of latter kind is the subject of this note.\\
	
	\noindent
	Under the assumption that $X$ is of non-negative Kodaira dimension and $k$ is of characteristic 0 it is proved in \cite[Theorem 4.2]{Fakhruddin-02} that there exists an abelian variety $A$ with a finite \'{e}tale map to $X$.
	Furthermore it was suggested in a footnote to Theorem 4.2 of \textit{loc. cit.} that a positive characteristic analogue should hold under certain additional hypothesis on the self map $\phi$. The following may be viewed as a result in this direction.

	\begin{theorem} \label{thm:S-is-0-dim-char-p}
		Let $X$ be a smooth projective variety of non-negative Kodaira dimension over an algebraically closed field $k$ of characteristic $p>0$, $L$ an ample line bundle on $X$ and $\phi : X \rightarrow X$ a separable map such that $\phi^*L \cong L^{\otimes d}$ for an integer $d>1$.
		Suppose that the \'{e}tale fundamental group of $X$ contains a normal abelian subgroup of finite index. Then there is an abelian variety $A$ with a finite \'{e}tale map to $X$.
	\end{theorem}
	
	\noindent
	The hypothesis on the fundamental group in the above theorem holds under the assumption that the variety $X$ admits a lift to characteristic 0 as shown in Lemma \ref{lem:virtually-abelian-fundamental-group}. Under the lifting assumption on $X$, the arguments of the proof of Theorem \ref{thm:S-is-0-dim-char-p} can be used to provide a proof of Theorem \ref{thm:S-is-0-dim-char-p} that uses results of more elementary nature (See Lemma \ref{lem:short-proof-thm-1}).
	Finally, we prove appropriate analogues of \cite[Theorem 4.2]{Fakhruddin-02} and Theorem \ref{thm:S-is-0-dim-char-p} over fields that are not necessarily algebraically closed.
	It is noteworthy that over characteristic 0 fields, the assumption on fundamental groups is known to be true for any smooth projective variety with polarized self map thanks to \cite[Corollary 4.6(2)]{Nakayama-Zhang-10} and \cite[Theorem 1.21]{Greb-Kebekus-Peternell-16}.\\
	
	\noindent
	The article is organized as follows : In \S \ref{subsec:prelim-lemma-proof-thm-1} we state two lemma that are required for the proof of Theorem \ref{thm:S-is-0-dim-char-p}. In \S \ref{subsubsec:proof-thm-1-with-Fp} we complete the proof of the theorem under the assumption that $k = \overline{\mathbb{F}}_p$ and in \S \ref{subsec:proof-of-thm-1-reduction-step} we reduce proving the theorem to the case of $k=\overline{\mathbb{F}}_p$. \S \ref{sec:two-lemma-on-hypothesis} consists of two Lemma whose contents have been described in the paragraph above. \S \ref{sec:non-algebraically-closed-field} generalizes (an appropriately modified version of) the above theorem to non-algebraically closed fields. At the end, we provide a remark on generalizations of the results of this note to ``int-amplified" self maps.

	\section{Proof of Theorem \ref{thm:S-is-0-dim-char-p}} \label{sec:proof-of-thm-1}
	\subsection{Preliminary lemmas} \label{subsec:prelim-lemma-proof-thm-1}
	For any map $\phi : X \rightarrow X$ we write $\phi^r$ to denote the composite of $\phi$ with itself $r$ times.
	
	%lemma
	\begin{lemma} \label{lemma:finite-etale-replacement}
		Let $X$ be a smooth projective variety over $k$ such that the \'{e}tale fundamental group of $X$ contains a normal abelian subgroup of finite index. Suppose $L$ is an ample line bundle on $X$ and $\phi : X \rightarrow X$ a finite \'{e}tale map such that $\phi^*L \simeq L^{\otimes d}$ for an integer $d>1$.
		There there exists a connected finite \'{e}tale cover $Z$ of $X$, an ample line bundle $M$ on $Z$ and a finite \'{e}tale map $\psi : Z \rightarrow Z$ such that $\psi^*M \simeq M^{\otimes e}$ for an integer $e>1$ and the fundamental group of $Z$ is abelian.
	\end{lemma}
	%proof
	\noindent
	\emph{Proof.}
	Fix a point $x_0 \in X$, and let $G \subset \pi^{\textup{\'{e}t}}_1(X, x_0)$ be a normal abelian subgroup of finite index. Let $Y \rightarrow X$ be the Galois \'{e}tale cover of $X$ given by the normal subgroup $G$ and $y_0 \in Y$ a point lying over $x_0$.
	In particular, $\pi^{\textup{\'{e}t}}_1(Y, y_0) \cong G$ is abelian. Consider the smooth projective (possibly disconnected) schemes $Y^{(r)}$ defined by the cartesian diagram
	%commutative diagram
	\[
	\begin{tikzcd}
		Y^{(r)} \arrow[r] \arrow[d] & Y \arrow[d] \\
		X \arrow[r, "\phi^r"] & X
	\end{tikzcd}
	\]
	The left downward arrow $Y^{(r)} \rightarrow X$ is finite \'{e}tale and of the same degree as $f = \deg(Y \rightarrow X)$ for all $r \geq 1$.
	Since $\pi^{\textup{\'{e}t}}_1(X, x_0)$ is topologically finitely generated\footnote{Using Lefschetz hyperplane theorem for fundamental groups it suffices to prove that $\pi^{\textup{\'{e}t}}_1(X, x_0)$ is topologically finite generated for a smooth projective curve $X$ over an algebraically closed field $k$.
		Any smooth projective curve admits a lifting to characteristic zero. Let $\widetilde{X}$ be a lift in characteristic 0. Using GAGA we conclude that $\pi^{\textup{\'{e}t}}_1(\widetilde{X}, \widetilde{x}_0)$ is topologically finitely generated.
		Finally surjectivity of specialization morphism for \'{e}tale fundamental group implies the same for $\pi^{\textup{\'{e}t}}_1(X, x_0)$.} there are only finitely many homomorphisms $\pi^{\textup{\'{e}t}}_1(X, x_0) \rightarrow H$ for any finite group $H$.
	Therefore there are only finite \'{e}tale covers of $X$ of degrees $\leq f$, in particular the collection of finite \'{e}tale covers of $X$ given by $\lbrace Y^{(1)}, Y^{(2)}, \dots \rbrace$ is finite.
	Hence there exists $n> m$ such that $Y^{(n)} \cong Y^{(m)}$ as \'{e}tale covers of $X$ and we get a cartesian diagram
	%commutative diagram
	\[
	\begin{tikzcd}
		Y^{(m)} \arrow[r, "\psi"] \arrow[d] & Y^{(m)} \arrow[d] \\
		X \arrow[r, "\phi^{n-m}"] & X
	\end{tikzcd}
	\]
	The map $\psi$ is surjective by construction and therefore the map $\pi_0(\psi) : \pi_0(Y^{(m)}) \rightarrow \pi_0(Y^{(m)})$ induced on connected components is bijective.
	After possibly replacing $\psi$ by $\psi^r$ for some integer $r$, we may assume that $\psi$ fixes each connected component of $Y^{(m)}$.
	Let $Z=$ a connected component of $Y^{(m)}$ and $M$ be the pullback of $L$ to $Z$ and $\psi : Z \rightarrow Z$ induced by the cartesian diagram above. We check that $\psi^*M \cong M^{e}$ where $e = d^{n-m}$.
	This completes the proof after noting that $\pi^{\textup{\'{e}t}}_1(Y)$ is abelian and $\pi^{\textup{\'{e}t}}_1(Z)$ is a subgroup of $\pi^{\textup{\'{e}t}}_1(Y)$ since $Z \rightarrow Y$ is a finite \'{e}tale map. \hfill \qed

	%lemma
	\begin{lemma}[{\cite[\S III, Lemma 5]{Katz-Lang-81}}] \label{lemma:ker-alb-finite}
		Let $X$ be a smooth projective variety over $k$ whose fundamental group is abelian. Let $alb_X : X \rightarrow Alb_X$ be the Albanese map. Then the kernel of the map $(alb_X)_*$ induced on fundamental groups is finite.
	\end{lemma}

	\subsection{Completing the proof of Theorem \ref{thm:S-is-0-dim-char-p}}
	We know from \cite[Lemma 4.1]{Fakhruddin-02} that $\phi$ is \'{e}tale. Without loss of generality we may replace $X$ by its finite \'{e}tale cover. Now using Lemma \ref{lemma:finite-etale-replacement} we may and do assume that the \'{e}tale fundamental group of $X$ is abelian.
	In light of \cite[Theorem 0.3(3)]{Hacon-Patakfalvi-Zhang-19} it suffices to prove that $alb_X$ is a finite morphism. This is what we will show.

	\subsubsection{Proof under the assumption $k = \overline{\mathbb{F}}_p$} \label{subsubsec:proof-thm-1-with-Fp}
	Since $k=\overline{\mathbb{F}}_p$, any point of $X(k)$ is preperiodic for (any self map) $\phi$ and hence there is a periodic point for $\phi$. Replacing $\phi$ by $\phi^n$ we may assume that there is a point $x_0 \in X$ such that $\phi(x_0) = x_0$. With $alb_X(x_0)$ as a base point, we will henceforth assume that $Alb_X$ is an abelian variety with $alb_X(x_0)$ as identity.
	It follows from universal property of the Albanese morphism that there exists a morphism $\phi_{Alb_X} : Alb_X \rightarrow Alb_X$ making the diagram below commute
	%commutative diagram
	\[
	\begin{tikzcd}
		X \arrow[r, "\phi"] \arrow[d, "alb_X"'] & X \arrow[d, "alb_X"]\\
		Alb_X \arrow[r, "\phi_{Alb_X}"'] & Alb_X
	\end{tikzcd}
	\]
	By construction $\phi_{Alb_X}(alb_X(x_0)) = alb_X(\phi(x_0)) = alb_X(x_0) = 0$, therefore $\phi_{Alb_X}(0) = 0$. Hence $\phi_{Alb_X}$ is a homomorphism of abelian varieties.
	Since $\im(alb_X)$ is not contained in any proper abelian subvariety, the commutative diagram above implies that $\im(\phi_{Alb_X})$ is not contained in a proper abelian subvariety. Thus $\phi_{Alb_X}$ is surjective homomorphism of abelian varieties of same dimension. That is $\phi_{Alb_X}$ is an isogeny, hence flat as a corollary.\\

	\noindent
	A result of Igusa (see for instance \cite[Lemma 1.3]{Mehta-Srinivas-87}) says that $alb^*_X : \mathrm{H}^0(Alb_X, \Omega^1_{Alb_X}) \rightarrow \mathrm{H}^0(X, \Omega^1_X)$ is injective. We have the following commutative diagram induced by the one above
	%commutative diagram
	\begin{equation} \label{eqn:commutative-diagram}
		\begin{tikzcd}
			\mathrm{H}^0(X, \Omega^1_X) & \arrow[l, "\phi^*"'] \mathrm{H}^0(X, \Omega^1_X)\\
			\mathrm{H}^0(Alb_X, \Omega^1_{Alb_X}) \arrow[u, "alb^*_X"] & \arrow[l, "\phi^*_{Alb_X}"]\mathrm{H}^0(Alb_X, \Omega^1_{Alb_X}) \arrow[u, "alb^*_X"']
		\end{tikzcd}
	\end{equation}
	Since $\phi$ is a surjective \'{e}tale(hence separable) map the top horizontal arrow is injective. The vertical arrows are injective as argued above. Therefore $\phi^*_{Alb_X}$ is injective, and hence isomorphism (since $\mathrm{H}^0(Alb_X, \Omega^1_{Alb_X})$ is finite dimensional).
	Since $\Omega^1_{Alb_X}$ is a free $\mathcal{O}_X$-module and $\Gamma(Alb_X, \mathcal{O}_{Alb_X}) = k$, it follows that the kernel and cokernel of the map $\phi^*_{Alb_X}\Omega^1_{Alb_X} \rightarrow \Omega^1_{Alb_X}$ are free $\mathcal{O}_X$-module. Applying $\phi_{Alb_X*}$ and using $R^i\phi_{Alb_X*} = 0$ for all $i>0$ we get an exact sequence below where $\Omega^1_{\phi_{Alb_X}}$ is the sheaf of relative differentials
	\[
	\phi_{Alb_X*}\phi^*_{Alb_X}\Omega^1_{Alb_X} \rightarrow \phi_{Alb_X*}\Omega^1_{Alb_X} \rightarrow \phi_{Alb_X*}\Omega^1_{\phi_{Alb_X}} \rightarrow 0.
	\]
	Note that the argument above implies that $\Omega^1_{\phi_{Alb_X}}$ is free. On the other hand the induced map on global section of the composite map $\Omega^1_{Alb_X} \rightarrow \phi_{Alb_X*}\phi^*_{Alb_X}\Omega^1_{Alb_X} \rightarrow \phi_{Alb_X*}\Omega^1_{Alb_X}$ is an isomorphism, hence $\Gamma(Alb_X, \phi_{Alb_X*}\Omega^1_{\phi_{Alb_X}}) = 0$. Since $\Omega^1_{\phi_{Alb_X}}$ is free (as observed above), $\dim_k\Gamma(Alb_X, \phi_{Alb_X*}\Omega^1_{\phi_{Alb_X}}) = $ rank of $\Omega^1_{\phi_{Alb_X}}$, therefore $\Omega^1_{\phi_{Alb_X}} = 0$.
	We have already seen that $\phi_{Alb_X}$ is flat, hence $\phi_{Alb_X}$ is finite \'{e}tale.\\
	
	\noindent
	Consider
	\[
	Z = \lbrace p \in Alb_X(k) \mid alb^{-1}_X(p) \text{ is not a finite set} \rbrace.
	\]
	It is clear that $Z$ is a Zariski closed subset, and we equip $Z$ with the reduced induced subscheme structure. Also the fact that $\phi, \phi_{Alb_X}$ is finite \'{e}tale and the diagram \eqref{eqn:commutative-diagram} is commutative implies that $\phi_{Alb_X}(Z) = Z$.
	Since $\phi_{Alb_X}$ is surjective after possibly replacing $\phi_{Alb_X}$ by $\phi^n_{Alb_X}$(which can be achieved by replacing $\phi$ by $\phi^n$) we may assume that $\phi_{Alb_X}$ preserves irreducible components of $Z$.
	Replacing $Z$ by any irreducible component of $Z$, we get that $Z$ is variety over $k$ and $\phi_{Alb_X}|_Z : Z \rightarrow Z$ is a finite surjective map ($\phi_{Alb_X}$ is finite).
	To prove that $alb_X$ is a finite map it suffices to prove that $Z = \emptyset$ since any proper quasi-finite map is finite. Suppose that $Z$ is nonempty. Then (arguing as in first paragraph of this subsection) there is a point $p \in Z(k)$ that is periodic for $\phi_{Alb_X}$ of order $n$ for some $n>0$. Replacing $\phi$ by $\phi^n$ (and as a consequence replacing $\phi_{Alb_X}$ by $\phi^n_{Alb_X}$) we may assume that $\phi_{Alb_X}(p) = p$.
	Let $Y$ be a connected component of $alb^{-1}_X(p)$, $y_0 \in Y$ be any point and $i : Y \rightarrow X$ be the inclusion map. We now make the following\\
	
	\noindent
	\textbf{Claim :} Let $i_* : \pi^{\textup{\'{e}t}}_1(Y, y_0) \rightarrow \pi^{\textup{\'{e}t}}_1(X, y_0)$ and $G = \im(i_*)$. Then $G$ is an infinite group.
	%proof
	\emph{Proof of claim.}
	Let $X_p = alb^{-1}_X(p)$. Since $\phi_{Alb_X}(p) = p$, $\phi_{Alb_X}$ is finite \'{e}tale and the diagram \eqref{eqn:commutative-diagram} is commutative we get that $\phi$ restricts to a finite \'{e}tale map
	\[
	\phi|_{X_p} : X_p \rightarrow X_p.
	\]
	Since there are finitely many connected components of $X_p$, we may replace $\phi$ by $\phi^r$ and assume that $\phi|_{X_p}$ fixes each connected component of $X_p$.
	Therefore for any connected component $Y$ of $X_p$ we get that $\phi|_Y : Y \rightarrow Y$ is a finite \'{e}tale map. It is also easy to see that $(\phi|_Y)^*L \simeq (L|_Y)^{\otimes d}$, and therefore $\phi|_Y$ is of degree $d^{\dim Y}$.
	Since the restriciton of $\phi^r$ to $Y$ provides finite \'e{tale} covers of $Y$ of degree $d^{r\dim Y}$ we get that $G = \im(i_*)$ must be infinite. \hfill \qed\\
	
	\noindent
	On the other hand $G$ lies in the kernel of the map $(alb_X)_* : \pi^{\textup{\'{e}t}}_1(X, x_0) \rightarrow \pi^{\textup{\'{e}t}}_1(Alb_X, 0)$ since $i\circ alb_X$ is constant. Lemma \ref{lemma:ker-alb-finite} implies that $G$ must be finite.
	Thus we get a contradiction to the assumption that $Z$ is nonempty. Thus completing the proof in this case.

	\subsubsection{Reduction to $k = \overline{\mathbb{F}}_p$} \label{subsec:proof-of-thm-1-reduction-step}
	Any $X, \phi$, and $L$ are defined over a finitely generated subfield of $k$, so we may assume that $k$ is of finite transcendence degree say $r$ over $\overline{\mathbb{F}}_p$.
	Note that $X,\phi, L$ and the isomorphism $\phi^*L \cong L^{\otimes d}$ is defined by finitely many equations. Hence there exists a d.v.r. whose field of fractions is $k$ and residue field is an algebraically closed of transcendence degree $r-1$ such that there is a triple $(\textbf{X}, \Phi, \textbf{L})$ consisting of smooth projective varieties $\textbf{X}$ over $\Spec(R)$, $\Phi : \textbf{X} \rightarrow \textbf{X}$, a relatively ample line bundle $\textbf{L}$ on $\textbf{X}$ such that $\Phi^*\textbf{L} \simeq \textbf{L}^{\otimes d}$ and over the generic point of $\Spec(R)$ the triple coincides with $(X, \phi, L)$.
	Let $(X_0, \phi_0, L_0)$ be the triple obtained by restricting the triple $(\textbf{X}, \Phi, \textbf{L})$ over $\Spec(R)$ to the closed point in $\Spec(R)$. Then, it follows from the construction that $(X_0, \phi_0, L_0)$ is a smooth projective polarized variety over an algebraically closed field of transcendence degree $r-1$. Note that the specialization map $\pi^{\textup{\'{e}t}}_1(X) \rightarrow \pi^{\textup{\'{e}t}}_1(X_0)$ is an isomorphism and hence $\pi^{\textup{\'{e}t}}_1(X_0)$ is abelian.\\
	
	\noindent
	Assume that Albanese map is known to be finite for all smooth projective polarized varieties $(X, \phi, L)$ defined over an algebraically closed of transcendence degree $r-1$ with $\phi$ finite \'{e}tale, and $\pi^{\textup{\'{e}t}}_1(X)$ is abelian.
	Consider the relative Albanese map $alb_{\textbf{X}} : \textbf{X} \rightarrow Alb_{\textbf{X}}$. Then $alb_{\textbf{X}}$ over the generic fiber coincides with $alb_X$ and over the special fiber coincides with $alb_{X_0}$.
	The induction hypothesis implies that $alb_{X_0}$ is a finite map, and hence so is $alb_X$. This completes the proof of the induction step.

	\section{Two lemma on the hypothesis of the Theorem} \label{sec:two-lemma-on-hypothesis}
	\begin{lemma} \label{lem:virtually-abelian-fundamental-group}
		Let $X$ be a smooth projective variety of non-negative Kodaira dimension over an algebraically closed field $k$, $L$ an ample line bundle on $X$ and $\phi : X \rightarrow X$ a separable map such that $\phi^*L \cong L^{\otimes d}$ for an integer $d>1$. If $X$ admits a lift to characteristic zero, then \'{e}tale fundamental group of $X$ contains a normal abelian subgroup of finite index.
	\end{lemma}
	
	\noindent
	\emph{Proof.}
	Let $\omega_X$ denote the canonical bundle of $X$. Let $\pi : \mathscr{X} \rightarrow \Spec(R)$ be a flat projective morphism, where $R$ is a discrete valuation ring with residue field $k$. Suppose $K$ is the fraction field of $R$. Since $\pi$ is flat and special fiber of $\pi$ is smooth, it follows that the morphism $\pi$ is smooth and proper.
	Let $j : \widetilde{X} \rightarrow \mathscr{X}$ be the geometric generic fiber of $\pi$ and $i : X \rightarrow \mathscr{X}$ be the special fiber of $\pi$.
	Using the functoriality of the chern class we get that $j^*c_1(\omega_{\mathscr{X}/\Spec(R)}) = c_1(\omega_{\widetilde{X}})$ and $i^*(c_1(\omega_{\mathscr{X}/\Spec(R)})) = c_1(\omega_X)$. Now \cite[Lemma 4.1]{Fakhruddin-02} implies that $\omega_X$ is torsion and therefore its first Chern class $c_1(\omega_X)$ in $\ell$-adic cohomology is 0. Since $\pi$ is a smooth proper morphism the maps $i^*$ and $j^*$ in $\ell$-adic \'{e}tale cohomology are isomorphisms. Therefore $c_1(\omega_{\widetilde{X}}) = 0$.
	Now, applying \cite[Th\'{e}or\`{e}me 2]{Beauville-83} we get that there exists an abelian variety $A$ and a simply connected variety $Z$ such that $A\times Z \rightarrow \widetilde{X}$ is a finite \'{e}tale cover. But $\pi^{\textup{\'{e}t}}_1(A\times Z) \cong \pi^{\textup{\'{e}t}}_1(A)$ is abelian.
	Therefore $H = \pi^{\textup{\'{e}t}}_1(A\times Z)$ is an abelian subgroup of $\pi^{\textup{\'{e}t}}_1(\widetilde{X})$ of finite index.
	Let $f_1, f_2, \dots, f_r \in \pi^{\textup{\'{e}t}}_1(\widetilde{X})$ be representatives of elements in $\pi^{\textup{\'{e}t}}_1(\widetilde{X})/H$.
	Then $G = \underset{i= 1}{\overset{r}{\cap}} f_iHf^{-1}_i$ is an abelian normal subgroup of $\pi^{\textup{\'{e}t}}_1(\widetilde{X})$ of finite index.
	We may conclude using the fact that the specialization map $\pi^{\textup{\'{e}t}}_1(\widetilde{X}) \rightarrow  \pi^{\textup{\'{e}t}}_1(X)$ is an isomorphism since $\pi : \mathscr{X} \rightarrow \Spec(R)$ is smooth. \hfill \qed

	\noindent
	Under the assumptions of the above Lemma, we provide a proof of Theorem \ref{thm:S-is-0-dim-char-p} that does not depend on the theorem of Hacon--Patakfalvi--Zhang.
	
	\begin{lemma} \label{lem:short-proof-thm-1}
		Let $(X, \phi, L)$ be as in the lemma above. Then there is an abelian variety with a finite \'{e}tale map to $X$. 	
	\end{lemma}
	
	\noindent
	\emph{Proof.}
	In light of Lemma \ref{lem:virtually-abelian-fundamental-group} and Lemma \ref{lemma:finite-etale-replacement} we may assume after replacing $X$ by a finite \'{e}tale cover that the fundamental group of $X$ is abelian and the triple $(X, \phi, L)$ continues to satisfy the assumptions of the theorem.
	Recall that $X$ admits a lift to characteristic zero; i.e. there is a smooth proper morphism $\pi : \mathscr{X} \rightarrow \Spec(R)$ where $R$ is a complete discrete valuation ring with $k$ as its residue field and $K$ as its field of fractions. Let $\widetilde{X}$ be generic fiber of $\pi$ and $\widetilde{a} : \widetilde{X} \rightarrow \widetilde{A}$ is the Albanese map.
	As an application of a result of Beauville (as in Lemma \ref{lem:virtually-abelian-fundamental-group}) we get that exists an abelian variety $\widetilde{B}$, a simply connected variety $Z$, a finite \'{e}tale map $\widetilde{h} : \widetilde{B} \times Z \rightarrow \widetilde{X}$.
	Moreover since $Alb_Z = 0$, the map $\widetilde{a}\circ \pi : \widetilde{B} \times Z \rightarrow \widetilde{A}$ factors through the projection $\widetilde{B} \times Z \rightarrow \widetilde{B}$. Therefore $im(\widetilde{a}\circ \widetilde{h}) \subset \widetilde{A}$ is an abelian subvariety.
	Since $im(\widetilde{a})$ is not contained in any proper abelian subvariety of $\widetilde{A}$ we get that $\im(\widetilde{a}\circ \widetilde{h}) = \widetilde{A}$. Hence $\dim(\widetilde{B}) \geq \dim(\widetilde{A})$.
	We get from the proof of the theorem that $\dim(X) = \dim(Alb_X)$. This equality implies that $\dim(\widetilde{X}) = \dim(X) = \dim(Alb_X) = \dim(\widetilde{A})$.
	Therefore the \'{e}tale cover $\widetilde{B} \times Z \rightarrow \widetilde{X}$ considered in previous paragraph satisfies $\dim(\widetilde{B}) + \dim(Z) = \dim(\widetilde{A})$. But $\dim(\widetilde{B}) \geq \dim(\widetilde{A})$ due to the equality $\im(\widetilde{a}\circ \widetilde{h}) = \widetilde{A}$, therefore we get $\dim(Z) = 0$ and $\dim(\widetilde{B}) = \dim(\widetilde{A})$. Hence $\widetilde{h} : \widetilde{B} \rightarrow \widetilde{X}$ is a finite \'{e}tale cover such that $\widetilde{a}\circ \widetilde{h}$ is a surjective map among abelian varieties with $\dim(\widetilde{B}) = \dim (\widetilde{A})$.
	Hence we may conclude that $\widetilde{a} : \widetilde{X} \rightarrow \widetilde{A}$ is a finite \'{e}tale map.
	Now Lang--Serre theorem implies that $\widetilde{X}$ is also an abelian variety. Let $\widetilde{0} \in \widetilde{X}(K)$ denote the identity element, and $e : \Spec(R) \rightarrow \mathscr{X}$ be the arrow obtained in the diagram below using valuative criterion of properness
	%commutative diagram
	\[
	\begin{tikzcd}
		\Spec(K) \arrow[r, "\widetilde{0}"] \arrow[d] & \mathscr{X} \arrow[d, "\pi"] \\
		\Spec(R) \arrow[ur, "e"] \arrow[r, equals] & \Spec(R)
	\end{tikzcd}
	\]
	Applying \cite[Theorem 6.14]{Mumford-GIT} we conclude that $\mathscr{X} \rightarrow \Spec(R)$ is an abelian scheme. Therefore its special fiber $X$ is an abelian variety. \hfill\qed

	\section{Results over non-algebraically closed fields} \label{sec:non-algebraically-closed-field}
	In this section we prove analogues of result \cite[Theorem 4.2]{Fakhruddin-02} and Theorem \ref{thm:S-is-0-dim-char-p} over fields that are not necessarily algebraically closed.
	
	\subsection{Preliminary lemmas}
	%lemma
	\begin{lemma} \label{lemma:descent-Beiberbach-group}
		Let $X$ be an algebraic variety over a field $k$, and let $\bar{x} \in X(\bar{k})$ be the geometric point lying over $x \in X(k)$ and $j_x$ be the splitting of the map $\pi^{\textup{\'{e}t}}_1(X, \bar{x}) \rightarrow \mathrm{Gal}(\overline{k}/k)$.
		Suppose that $\pi^{\textup{\'{e}t}}_1(X_{\overline{k}}, \bar{x})$ has a characteristic subgroup $H$ of finite index. Let $Y \rightarrow X$ be the \'{e}tale cover given by the subgroup $H\cdot j_x(\mathrm{Gal}(\overline{k}/k))$, then the \'{e}tale cover $Y_{\overline{k}} \rightarrow X_{\overline{k}}$ obtained via base change is given by the subgroup $H \subset \pi^{\textup{\'{e}t}}(X_{\overline{k}}, \bar{x})$.
	\end{lemma}
	\noindent
	\emph{Proof.}
	The Galois group $\mathrm{Gal}(\overline{k}/k)$ acts on the algebraic fundamental group $\pi^{\textup{\'{e}t}}_1(X_{\overline{k}}, \bar{x})$ by outer automorphisms.
	But any automorphism preserves the characteristic subgroup and hence it preserves the subgroup $H$. Thus the group generated by $H$ and $j_x(\mathrm{Gal}(\overline{k}/k))$ is a normal subgroup of $\pi^{\textup{\'{e}t}}_1(X_{\overline{k}}, \bar{x})$.
	Let $\pi : Y \rightarrow X$ be an \'{e}tale cover given by $H\cdot j_x(\mathrm{Gal}(\overline{k}/k))$. Using the functoriality of the fundamental exact sequence of \'{e}tale fundamental groups it is easy to see  that the map induced by $\pi$ on algebraic fundamental group is the inclusion $H \subset \pi^{\textup{\'{e}t}}(X_{\overline{k}}, \bar{x})$.
	This proves the last assertion of the lemma.\hfill\qed
	
	\noindent
	The following lemma is a corollary of \cite[Theorem 1]{Khukhro-Makarenka-07}.
	
	\begin{lemma} \label{lem:virtually-abelian-characteristic}
		Any virtually abelian group contains a characteristic abelian subgroup of finite index.
	\end{lemma}

	\subsection{Char(k)=0}
	
	\begin{theorem}
		Let $X$ be a smooth projective variety of non-negative Kodaira dimension over a field $k$ of characteristic 0, $L$ an ample line bundle on $X$ and $\phi : X \rightarrow X$ a self map of $X$ such that $\phi^*L \cong L^{\otimes d}$ for an integer $d>1$. Assume that $X(k) \neq \emptyset$, then there is an abelian variety $A$ with a finite \'{e}tale map to $X$ over $k$.
	\end{theorem}
	
	\noindent
	It is clear that the assumption of $X(k) \neq \emptyset$ is required on account of existence of para-abelian varieties.\\

	\noindent
	\emph{Proof}
	Let $x \in X(k)$ and $\overline{x} \in X(\overline{k})$ be a geometric point lying over $x$. Recall from \cite[Theorem 4.2]{Fakhruddin-02} that $X_{\overline{k}}$ admits a finite \'{e}tale cover by an abelian variety $B$. Thus the image of algebraic fundamental group of $B$ in $\pi^{\textup{\'{e}t}}_1(X_{\overline{k}}, \overline{x})$ say $N$ is an abelian normal subgroup. In particular $\pi^{\textup{\'{e}t}}_1(X_{\overline{k}}, \bar{x})$ is virtually abelian and hence by Lemma \ref{lem:virtually-abelian-characteristic} contains a characteristic abelian subgroup $\Gamma$ of finite index.
	Moreover, it can be arranged that $\Gamma \subset N$, by simply replacing $\Gamma$ by $h\Gamma \subset \Gamma$ for an integer $h = |\Gamma/(N\cap \Gamma)|$.
	Consider the \'{e}tale cover $Z \rightarrow X_{\overline{k}}$ given by the subgroup $\Gamma$ of $\pi^{\textup{\'{e}t}}_1(X_{\overline{k}}, \overline{x})$ and we denote the quotient $\pi^{\textup{\'{e}t}}_1(X_{\overline{k}}, \overline{x})/\Gamma$ by $G$.
	Since $\Gamma \subset N$, hence $Z$ is a a finite \'{e}tale cover of $B$ and hence itself is an abelian variety by Lang--Serre theorem.
	We may apply Lemma \ref{lemma:descent-Beiberbach-group} to get that there exists a smooth projective variety $A$ defined over $k$ such that $A \rightarrow X$ is a $G$-torsor and $A_{\overline{k}}$ is the \'{e}tale cover of $X_{\overline{k}}$ given by the subgroup $\Gamma$ which is known to be an abelian variety from the discussion above.
	To complete the proof it remains to show that there exists a form $A'$ of $A$ (over $k$) equipped with a finite \'{e}tale map $A' \rightarrow X$, and $A'(k) \neq \emptyset$.
	
	\noindent
	We recall a few relevant facts from \cite{Harari-Skorobogatov-02}. We say that $\alpha \in \mathrm{Aut}(A_{\overline{k}}/k)$ is a $k$-semilinear automorphism of $A_{\overline{k}}$ if for some $g \in \mathrm{Gal}(\overline{k}/k)$ the diagram below commutes
	\begin{equation} \label{eqn:diagram}
		\begin{tikzcd} A_{\overline{k}} \arrow[r, "\alpha"] \arrow[d] & A_{\overline{k}} \arrow[d] \\ \Spec(\overline{k}) \arrow[r, "(g^*)^{-1}"] & \Spec(\overline{k}) \end{tikzcd}.
	\end{equation}
	Denote the subset of $k$-semilinear automorphism of $A_{\overline{k}}$ to be $S\mathrm{Aut}(A_{\overline{k}})$. The map which sends $\alpha \in S\mathrm{Aut}(A_{\overline{k}})$ to $g \in \mathrm{Gal}(\overline{k}/k)$ where $\alpha$ and $g$ are as in the diagram above defines a homomorphism
	\[
	q : S\mathrm{Aut}(A_{\overline{k}}/X) \rightarrow \mathrm{Gal}(\overline{k}/k).
	\]
	Let $E = S\mathrm{Aut}(A_{\overline{k}}/X) \colonequals \mathrm{Aut}(A_{\overline{k}}/X) \cap S\mathrm{Aut}(A_{\overline{k}}/k)$.
	Since $A$ is defined over $k$ the map $q : E \rightarrow \mathrm{Gal}(\overline{k}/k)$ is split surjective. We have the following short exact sequence
	\[
	1 \rightarrow \mathrm{Aut}(A_{\overline{k}}/X_{\overline{k}}) \rightarrow E \xrightarrow{q} \mathrm{Gal}(\overline{k}/k) \rightarrow 1.
	\]
	Applying \cite[Theorem 2.5]{Harari-Skorobogatov-02} with the choice of $e \in A(\overline{k})$ being the identity element, we get that there exists a (homomorphic) section $j_{e} : \mathrm{Gal}(\overline{k}/k) \rightarrow E$ of the map $q : E \rightarrow \mathrm{Gal}(\overline{k}/k)$ such that $e$ is fixed under $j_e(\mathrm{Gal}(\overline{k}/k))$.
	Consider $A'$(resp. $G'$) the form of $A_{\overline{k}}$(resp. $G_{\overline{k}}$) defined by $j_{e}$ (see \cite[Lemma 1.4]{Harari-Skorobogatov-02}). By construction $A'(k) \neq \emptyset$ and the map $A_{\overline{k}} \rightarrow X$ factors through $A'$ to give a morphism $A' \rightarrow X$.
	Theorem 2.5 of \textit{loc. cit.} moreover implies that the map $A' \rightarrow X$ is a torsor under the finite \'{e}tale group scheme $G'$. This completes the proof. \hfill\qed

	\subsection{Char(k)$>$0}
	\begin{theorem}
		Let $X$ be a smooth projective variety of non-negative Kodaira dimension over a field $k$ of positive characteristic, $L$ an ample line bundle on $X$ and $\phi : X \rightarrow X$ a separable map such that $\phi^*L \cong L^{\otimes d}$ for an integer $d>1$. Assume that the algebraic fundamental group of $X$ is virtually abelian and $X(k) \neq \emptyset$, then there is an abelian variety $A$ with a finite \'{e}tale map to $X$ over $k$.
	\end{theorem}

	\noindent
	The proof in this case follows verbatim as in the case of characteristic zero by applying Theorem \ref{thm:S-is-0-dim-char-p} instead of \cite[Theorem 4.2]{Fakhruddin-02}.\\

	\noindent
	\emph{Remark.}
	Let $X$ be a smooth projective variety over an algebraically closed field $k$. A morphism $f : X \rightarrow X$ is said to be int-amplified if there exists an ample line bundle $L$ such that $f^*L\otimes L^{-1}$ is an ample line bundle.
	It is easy to see that $f$ is a finite map and hence is surjective in this case. Arguing as before, we may conclude that $f$ is \'{e}tale. Finally note that $f$ is not an automorphism and hence $\deg(f)>1$. The proof of Theorem \ref{thm:S-is-0-dim-char-p} can be adapted using the observations made above to extend the theorem to the case of smooth projective varieties with int-amplified self map $f$.
	The argument of this note can be extended to the int-amplified case if the triviality of $c_1(\omega_{X})$ can be shown.

\section*{Acknowledgements}
The author thanks N. Fakhruddin and V. Srinivas for helpful discussions, comments, and suggestions.

\end{document}